\documentclass[11pt]{article}
\usepackage{amsmath,amssymb,amsfonts,amsthm,graphicx}
\usepackage{hyperref}
\newtheorem{theorem}{Theorem}[section]
\newtheorem{proposition}[theorem]{Proposition}
\newtheorem{corollary}[theorem]{Corollary}
\newtheorem{lemma}[theorem]{Lemma}
\newtheorem{remark}[theorem]{Remark}

\numberwithin{equation}{section}
\title{A two-dimensional Gauss-Kuzmin theorem for $N$-continued fraction expansions}
\author{
    Gabriela Ileana Sebe\footnote{e-mail: igsebe@yahoo.com.} \\
    \emph{\small Politehnica University of Bucharest, Faculty of Applied Sciences},\\
    \emph{\small Splaiul Independentei 313, 060042, Bucharest, Romania} and \\
    \emph{\small Institute of Mathematical Statistics and Applied Mathematics}, \\
     \emph{\small Calea 13 Sept. 13, 050711 Bucharest, Romania} \\
    and\\
    Dan Lascu\footnote{e-mail: lascudan@gmail.com.}\nonumber \\
    \emph{\small Mircea cel Batran Naval Academy, 1 Fulgerului, 900218 Constanta,
    Romania} \\
    }
\begin{document}
\maketitle
\thispagestyle{empty}
\begin{abstract}

A two-dimensional Gauss-Kuzmin theorem for $N$-continued fraction expansions is shown. More exactly, we obtain a Gauss-Kuzmin theorem related to the natural extension of the measure-dynamical system corresponding to these expansions.
Then, using characteristic properties of the transition operator associated with the random system with complete connections underlying $N$-continued fractions on the Banach space of complex-valued functions of bounded variation we derive explicit lower and upper bounds for the convergence rate of the distribution function to its limit.
\end{abstract}

\section{Introduction}

The purpose of this paper is to show a two-dimensional Gauss-Kuzmin theorem for $N$-continued fraction expansions introduced by Burger et al. \cite{Burger-2008}, and studied by Dajani et al. \cite{DKW-2013}.

In the last twenty-five years some versions of the so-called \textit{two-dimensional theorem} related to the \textit{Gauss-Kuzmin problem} are obtained for different types of continued fraction expansions.
For example, the case of the regular continued fraction (RCF) expansion was extensively studied in \cite{DK-1994}, \cite{I-1997} and \cite{IK-2002}.
For Hurwitz' singular continued fractions there are known the results obtained in \cite{DK-1999} and \cite{Sebe-2000}. Also, in \cite{Sebe-2001} Sebe proved the first two-dimensional Gauss-Kuzmin theorem which leads to an estimate of the approximation error by the expansion algorithm in the grotesque continued fraction.

Our study concludes the series of papers \cite{Lascu-2016} and \cite{Lascu-2017} where the second author discussed some metric properties of these expansions. Moreover he investigated the associated Perron-Frobenius operator and proved a Gauss-Kuzmin theorem for $N$-continued fraction applying the method of random systems with complete connections (RSCC) by Iosifescu \cite{IG-2009}.

The outline of this paper is as follows.
In Section 2 we present notions and preliminary results to be used in the sequel. We mention that the infinite-order-chain representation of the sequence of the partial quotients of the $N$-continued fraction expansion allows a concise formulation of the results obtained.
Section 3 presents a Gauss-Kuzmin theorem related to the natural extension \cite{Nakada-1981} of the measure-dynamical system corresponding to $N$-continued fraction expansions.
Section 4 is devoted to derive explicit lower and upper bounds of the error which provide a more refined estimate of the convergence rate. The key role in this section is played by the transition operator associated with the RSCC underlying $N$-continued fraction on the Banach space of function of bounded variation.
The last section collects some concluding remarks and give interesting numerical calculations.

\section{Preliminaries}

Fix an integer $N \geq 1$. The measure-theoretical dynamical system $(I,{\mathcal B}_{I},T_N)$ is defined as follows:
$I:=[0,1]$, $\mathcal{B}_I$ denotes the $\sigma$-algebra of all Borel subsets of $I$,
and $T_N$ is the transformation
\begin{equation}
T_N:I \to I; \quad
T_{N}(x):=
\left\{
\begin{array}{ll}
{\displaystyle \frac{N}{x}- \left\lfloor\frac{N}{x}\right\rfloor}&
{ \mbox{if } x \in I,}\\
0& \mbox{if }x=0.
\end{array}
\right. \label{1.1}
\end{equation}
In addition, we write $(I,{\mathcal B}_{I},G_N,T_N)$ as $(I,{\mathcal B}_{I}, T_N)$
with the following probability measure $G_N$ on $(I,{\mathcal B}_{I})$:
\begin{equation}
G_N (A) :=
\frac{1}{\log \frac{N+1}{N}} \int_{A} \frac{\mathrm{d}x}{x+N},
\quad A \in {\mathcal{B}}_I. \label{1.2}
\end{equation}
The measure $G_N$ is $T_N$-invariant, i.e., $G_N\left(T_N^{-1}(A)\right)=G_N(A)$ for any $A \in {\mathcal{B}}_I$ and the dynamical system
$(I,{\mathcal B}_{I},G_N,T_N)$ is ergodic.
For any $0<x<1$ put $a_n(x)=a_1\left(T_N^{n-1}(x)\right)$, $n \in \mathbb{N}_+:=\{1,2, \ldots\}$, with $T_N^{0}(x)=x$.
Then every irrational $0<x<1$ can be written in the form
\begin{equation}
x = \displaystyle \frac{N}{a_1+\displaystyle \frac{N}{a_2+\displaystyle \frac{N}{a_3+ \ddots}}} :=[a_1, a_2, a_3, \ldots]_N \label{1.3}
\end{equation}
where $a_n$'s are non-negative integers, $a_n \geq N$, $n \in \mathbb{N}_+$. We will call (\ref{1.3}) the \textit{$N$-continued fraction expansion} of $x$ and
$p_n(x) / q_n(x) := [a_1, a_2, a_3, \ldots, a_n]_N$ the \textit{$n$-th order convergent} of $x \in I$. Then $p_n(x) / q_n(x) \rightarrow x$, $n \rightarrow \infty$.
Here $p_n$'s and $q_n$'s satisfy for $n \in \mathbb{N}_+$ the following:
\begin{eqnarray}
p_n(x) &:=& a_n p_{n-1} + N p_{n-2}, \quad n \geq 2 \label{1.4} \\
q_n(x) &:=& a_n q_{n-1} + N q_{n-2}, \quad n \geq 1 \label{1.5}
\end{eqnarray}
with $p_0(x):=0$, $q_0(x) := 1$, $p_{-1}(x):=1$, $q_{-1}(x):=0$, $p_1(x):=N$, $q_1(x):=a_1$.
One easily shows that
\begin{equation}
\left| x - \frac{p_n(x)}{q_n(x)} \right| < \frac{N^n}{q^2_n(x)} < \frac{1}{N^n}, \quad n \geq 1. \label{1.6}
\end{equation}
For any $n \in \mathbb{N}_+$ and $i^{(n)}=(i_1, \ldots, i_n) \in \Lambda^n$, $\Lambda := \{N, N+1, \ldots \}$,
we will say that
\begin{equation}
I_N \left(i^{(n)}\right) = \{x \in I \setminus {\mathbb{Q}}:  a_k(x) = i_k \mbox{ for } k=1,\ldots,n \} \label{1.7}
\end{equation}
is the {\it fundamental interval of rank $n$} and make the convention that $I_N (i^{(0)}) = I \setminus {\mathbb{Q}}$.
For example, for any $i \in \Lambda$, we have
\begin{equation}
I_N(i) \,= \,\{x \in I \setminus {\mathbb{Q}}: a_1 = i \}
\,=\,\left( I \setminus {\mathbb{Q}}\right) \cap \left(\frac{N}{i+1}, \frac{N}{i}\right).\label{1.7'}
\end{equation}
For any $n \in \mathbb{N}_+$ we have
\begin{equation}
\lambda (T^n_N < x |a_1,\ldots, a_n )
= \frac{(s_n + N)x}{s_n x+N}, \quad x \in I \label{1.8}
\end{equation}
where $s_n := N \frac{q_{n-1}}{q_{n}}=[a_n, a_{n-1}, \ldots, a_2, a_1]_N$, $n \geq 1$, $s_0:=0$ and $\lambda$ is the Lebesgue measure on $I$.
Equation (\ref{1.8}) is the Brod\'en-Borel-L\'evy formula for this type of expansion. It allows us to determine the probability distribution of $(a_n)_{n \in \mathbb{N}_+}$ under $\lambda$. For any $i \in \Lambda$, we have
\begin{equation}
\lambda(a_1=i) = \frac{N}{i(i+1)}, \quad
\lambda\left(a_{n+1}=i |a_1,\ldots, a_n \right) = V_{N,i}(s_n) \label{1.9}
\end{equation}
where
\begin{equation}
V_{N,i}(x) := \frac{x+N}{(x+i)\,(x+i+1)}. \label{1.10}
\end{equation}
Hence the invariant probability measure of the transformation $T_N$ is $G_N$, the sequence $(a_n)_{n \in \mathbb{N}_+}$ is strictly stationary on $(I,{\mathcal B}_{I},G_N)$.
As such, a doubly infinite version of it, say $(\overline{a}_l)_{l \in \mathbb{Z}}$, should exist on a richer probability space.
Indeed, such a version can be effectively constructed on the natural extension $(I^2, {\mathcal B}_{I^2}, \overline{G_N}, \overline{T_N})$ of $(I,{\mathcal B}_I, G_N, T_N)$ where the transformation $\overline{T_N}$ is defined on the square space $(I^2,{\mathcal B}^2_I):=(I,{\mathcal B}_I) \times (I,{\mathcal B}_I)$ by
\begin{equation} \label{1.11}
\overline{T_N}(x,y) := \left( T_N(x), \, \frac{N}{a_1(x)+y} \right)
\end{equation}
and the extended measure $\overline{G_N}$ is defined by
\begin{equation} \label{1.12}
\overline{G_N}(B) :=\frac{1}{ \log\left(\frac{N+1}{N}\right)} \int\!\!\!\int_{B}
\frac{N\mathrm{d}x\mathrm{d}y}{(xy+N)^2}, \quad B \in {\mathcal{B}}^2_I.
\end{equation}
The transformation $\overline{T_N}$ is bijective on $I^2$ with the inverse
\begin{equation} \label{1.13}
(\overline{T_N})^{-1}(x, y) = \left(\frac{N}{a_1(y)+x}, \, T_N(y)\right), \quad (x, y) \in I^2.
\end{equation}
Iterations of (\ref{1.11}) and (\ref{1.13}) are given as follows for each $n \geq 2$:
\begin{eqnarray*}
\left(\overline{T_N}\right)^n(x, y) &=&
\left(T^n_N(x), [a_n(x), a_{n-1}(x), \ldots, a_2(x), a_1(x)+ y ]_N \right) \\
\left(\overline{T_N}\right)^{-n}(x, y) &=&
\left([a_n(y), a_{n-1}(y), \ldots, a_2(y), a_1(y)+x ]_N, T^{n}_N(y) \right).
\end{eqnarray*}
We define {\it extended incomplete quotients} $\overline{a}_l(x,y)$, $l \in \mathbb{Z}$ at $(x, y) \in I^2$ by
\begin{equation} \label{1.14}
\overline{a}_{l}(x, y) := \overline{a}_1\left(\,(\overline{T_N})^{l-1} (x, y) \,\right),
\quad l \in \mathbb{Z},
\end{equation}
with $\overline{a}_{1}(x, y)=a_1(x)$.
Clearly, for $n \in \mathbb{N}_+$ and any $(x, y) \in I^2$, $\overline{a}_n(x, y)=a_n(x)$, $\overline{a}_0(x, y)=a_1(y)$, $\overline{a}_{-n}(x, y)=a_{n+1}(y)$.
Since the extended measure $\overline{G_N}$ is preserved by $\overline{T_N}$, that is, $\overline{G_N}=\overline{G_N}(\overline{T_N})^{-1}$, the doubly infinite sequence $(\overline{a}_l(x,y))_{l \in \mathbb{Z}}$
is strictly stationary under $\overline{G_N}$.
The dependence structure of $(\overline{a}_l)_{l \in \mathbb{Z}}$ follows from the fact that
\begin{equation} \label{1.15}
\overline{G_N} ( [0, x] \times I \,|
\,\overline{a}_0, \overline{a}_{-1}, \ldots )
= \frac{(N+a)x}{ax + N} \quad \overline{G_N} \mbox{-}\mathrm{a.s.}
\end{equation}
for any $x \in I$, where $a:= [\overline{a}_0, \overline{a}_{-1}, \ldots]_N$. Hence
\begin{equation} \label{1.16}
\overline{G_N} (\left.\overline{a}_1 = i\right| \overline{a}_0, \overline{a}_{-1}, \ldots) = V_{N,i}(a) \quad \overline{G_N} \mbox{-}\mathrm{a.s.}
\end{equation}
where $V_{N,i}$ is as in (\ref{1.10}). We thus see that $\left(\overline{a}_l\right)_{l \in \mathbb{Z}}$ is an infinite-order-chain in the theory of dependence with complete connections \cite{IG-2009}.

Put $\overline{s}_\ell = [\overline{a}_{\ell}, \overline{a}_{\ell-1},\ldots]_N$, $\ell \in \mathbb{Z}$.
Note that $\overline{s}_\ell=N / (\overline{s}_{\ell-1}+\overline{a}_{\ell})$, $\ell \in \mathbb{Z}$.
It follows from (\ref{1.16}) that $\left( \overline{s}_{\ell} \right)_{\ell \in \mathbb{Z}}$ is a strictly stationary $I$-valued Markov process on $(I^2, {\mathcal{B}}^2_I, \overline{G_{N}})$ with the following transition mechanism:
from state $\overline{s} \in I$ the only possible transitions are those to states $N/(\overline{s}+i)$, $i \in \Lambda$, the transition probability being $V_{N,i}(\overline{s})$.
Clearly, whatever $\ell \in \mathbb{Z}$ we have $\overline{G_{N}}\left(\overline{s}_{\ell}<x\right)=\overline{G_{N}}\left(\overline{s}_{0}<x\right)=\overline{G_{N}}(I \times [0,x])=G_{N}([0,x])$, $x \in I$.

Motivated by (\ref{1.15}) we shall consider the family of (conditional) probability measures $(G_{N,a})_{a}$ on ${\mathcal{B}}_I$ defined by their distribution functions
\begin{equation} \label{1.17}
G_{N,a}([0,x]) = \frac{(N+a)x}{ax + N}, \quad x \in I, \, a \in I.
\end{equation}
In particular, $G_{N,0}$ is the Lebesgue measure on $I$.
For any $a \in I$ put $s_{0,a}=a$ and $s_{n,a} = N / \left(s_{n-1,a}+a_n\right)$, $n \in \mathbb{N}_+$.
These facts lead us to the random system with compete connections
$\{(I, {\mathcal{B}}_I), (\Lambda, \mathcal{P}(\Lambda)), u, V\}$,
where $\mathcal{P}(\Lambda)$ is the power set of $\Lambda$, $u:I \times \Lambda \rightarrow I$ is defined as
\begin{equation} \label{1.17'}
u(x,i):=u_i(x)=u_{N,i}(x)=\displaystyle \frac{N}{x+i}
\end{equation}
and
$V:I \times \Lambda \rightarrow I$ is defined as
\begin{equation} \label{1.17''}
V(x,i):= V_{N,i}(x)=v_i(x)
\end{equation}
with $V_{N,i}$ as in (\ref{1.10}), for all $x \in I$ and $i \in \Lambda$.

Then $\left(s_{n,a}\right)_{n \in \mathbb{N}_+}$ is an $I$ - valued Markov chain on $(I, {\mathcal{B}}_I, G_{N,a})$ which starts from $s_{a,0} = a$, $a \in I$, and has the following transition mechanism: from state $s \in I$ the only possible transitions are those to states $N/(s+i)$ with the corresponding transition probability $V_{N,i}(s)$, $i \in \Lambda$.
Let $B(I)$ denote the Banach space of all bounded $I$-measurable complex-valued functions defined on $I$ which is a Banach space under the supremum norm.
The transition operator of $(s_{n,a})_{n \in \mathbb{N}_+}$ takes $f \in B(I)$ into the function defined by
\begin{equation*}
E_a\left( \left. f(s_{n+1,a})\right| s_{n,a} = s \right) = \sum_{i \in \Lambda}V_{N,i}(s)f\left(\frac{N}{s+i}\right) = \{U f\}(s) \quad \mbox{for any } s \in I,
\end{equation*}
where $E_a$ stands for the mean-value operator with respect to the probability measure $G_{N,a}$, whatever $a \in I$, and $U$ is the Perron-Frobenius operator of $(I,{\mathcal B}_{I},G_N,T_N)$ \cite{Lascu-2017}.

In connection with the operator $U$, if we define $U^{\infty} f = \int_{I} f(x) G_N(\mathrm{d}x)$, $f \in B(I)$,
then we have
\begin{equation}
U^{\infty} U^{n} f = U^{\infty} f, \mbox{ for any } f \in B(I) \mbox{ and } n \in \mathbb{N}_+. \label{1.18}
\end{equation}
Note that
\begin{equation*}
G_{N,a}\left(A| a_1, \ldots, a_m\right) = G_{N,s_{m,a}}\left(T^{m}_N(A)\right),
\end{equation*}
for all $a \in I$, $A \in \sigma(a_{m+1}, \ldots)$ and $m \in \mathbb{N}_+$. In particular, it follows that the Brod\'en-Borel-L\'evy formula holds under $G_{N,a}$ for any $a \in I$, that is,
\begin{equation}
G_{N,a} (T^m_N < x |a_1,\ldots, a_m )
= \frac{(s_{m,a} + N)x}{s_{m,a} x+N}, \quad x \in I, m \in \mathbb{N}_+. \label{1.19}
\end{equation}
Note also that $(s_{n,a})_{n \in \mathbb{N}}$ under $G_{N,a}$ is a version of $(\overline{s}_{n})_{n \in \mathbb{N}}$ under $\overline{G_{N}}(\cdot | \overline{s}_0=a)$ for any $a \in I$.
\section{A two-dimensional Gauss-Kuzmin theorem}

The problem of finding the asymptotic behaviour of $T_N^{-n}(A)$ as $n \rightarrow \infty$, $A \in {\mathcal B}_{I}$, represents the Gauss-Kuzmin-type problem for $N$-continued fraction expansions. We define the functions $(F_n)_{n \in \mathbb{N}}$ on $I$ by
\begin{equation}
F_0(x):=x, \quad F_n(x) := \lambda(T^n_N < x), n \geq 1.  \label{2.1}
\end{equation}
The essential argument of the proof is the Gauss-Kuzmin-type equation which in this case is
\begin{equation}
F_{n+1}(x) = \sum_{i \in \Lambda}\left(F_n\left(\frac{N}{i}\right) - F_n\left(\frac{N}{i+1}\right) \right) \label{2.2}
\end{equation}
for $x \in I$ and $n \in \mathbb{N}$. The measure $G_N$ defined in (\ref{1.2}) is an eigenfunction of (\ref{2.2}), namely, if we put $F_n(x) = \log \left(\frac{x+N}{N}\right)$, $x \in I$, we obtain $F_{n+1}(x) = \log \left(\frac{x+N}{N}\right)$.
Let us recall now a slightly modified version of Gauss-Kuzmin theorem for $T_N$ proved in \cite{Lascu-2016}.
\begin{theorem} \label{G-K1}
Let $(I,{\mathcal B}_I,G_N, T_N)$ and $F_n$ be as above.
Then for any $x \in I$,
\begin{equation} \label{2.3}
F_n(x) = G_N([0,x]) + {\mathcal O} \left(q^n\right)
\end{equation}
with $0<q<1$.
\end{theorem}

In this section we will derive a Gauss-Kuzmin theorem related to the natural extension $(I^2, {\mathcal B}^2_{I}, \overline{G_N}, \overline{T_N})$ defined in Section 2.

For any $n \in \mathbb{N}_+$ and $x,y \in I$, let us define $\Delta_{x,y}=[0,x] \times [0,y]$ and the functions $F_n(x,y)$ by
\begin{equation} \label{2.4}
F_n(x,y) := \overline{\lambda} \left( \left(\overline{T_N}\right)^n \in \Delta_{x,y} \right),
\end{equation}
where $\overline{\lambda}$ is the Lebesgue measure on $I^2$.
Then the following holds.
\begin{theorem} \label{G-K2}
(A Gauss-Kuzmin theorem for $\overline{T_N}$)
For every $n\geq2$ and $(x,y) \in I^2$ one has
\begin{equation}\label{2.5}
F_n(x,y) = \frac{1}{\log\left(\frac{N+1}{N}\right)} \log \left( \frac{xy+N}{N}\right) + {\mathcal O} \left(q^n\right)
\end{equation}
with $0<q<1$.
\end{theorem}
For any $0 < y \leq 1$, put $\ell_1 := \left\lfloor\frac{N}{y}\right\rfloor$. Then $\left(\overline{T_N}\right)^{n+1}(x,y) \in \Delta_{x,y}$ is equivalent to
\begin{equation*}
\left(\overline{T_N}\right)^{n} \in \left( \bigcup_{i \geq \ell_1+N} \left[ \frac{N}{x+i}, \frac{N}{i} \right] \times [0,1] \right)
\cup \left( \left[\frac{N}{x+\ell_1}, \frac{N}{\ell_1} \right] \times \left[ \frac{N}{y}-\ell_1, 1 \right] \right).
\end{equation*}
Now, from (\ref{2.4}) we have the following recursion formula:
\begin{eqnarray}
  F_{n+1}(x,y) &=& \sum_{i \geq \ell_1} \left( F_n\left(\frac{N}{i},1\right) - F_n\left(\frac{N}{x+i},1\right) \right) \nonumber \\
               &-& \left( F_n\left(\frac{N}{\ell_1},\frac{N}{y}-\ell_1\right) - F_n\left(\frac{N}{x+\ell_1},\frac{N}{y}-\ell_1\right) \right). \label{2.7}
\end{eqnarray}
The measure $\overline{G_N}$ defined in (\ref{1.12}) is an eigenfunction of (\ref{2.7}), namely,
if we put $F_n(x,y) = \log \left(\frac{xy+N}{N}\right)$, $x,y \in I$, we obtain $F_{n+1}(x,y) = \log \left(\frac{xy+N}{N}\right)$.
\begin{lemma} \label{lema2.3}
Let $n \in \mathbb{N}$, $n \geq 2$ and let $y \in I \cap \mathbb{Q}$ with $y = [\ell_1, \ldots, \ell_d]_N$, $\ell_1, \ldots, \ell_d \in \Lambda$, $\ell_d \geq N+1$, where $d \leq \lfloor n/(N+1) \rfloor$. Then for every $x, x^* \in [0,1)$ with $x^* < x$,
\begin{equation*}
\left| F_n(x,y) - F_n(x^*,y) - \frac{1}{\log\left(\frac{N+1}{N}\right)} \log \left( \frac{xy+N}{x^*y+N} \right) \right| <
C \overline{\lambda}(\Delta_{x,y} \setminus \Delta_{x^*,y}) q^{n-d},
\end{equation*}
where $C$ is an universal constant.
\end{lemma}
\noindent \textbf{Proof.}
Let $y_0=y$, $y_i:=[\ell_{i+1}, \ldots, \ell_1]_N$, $i=1, \ldots, d-1$, and $y_d=0$.
Applying (\ref{2.7}) one gets
\begin{eqnarray}
F_{n}(x,y) - F_{n}(x^*,y) &=& \sum_{i \geq \ell_1} \left( F_{n-1}\left(\frac{N}{x^*+i},1\right) - F_{n-1}\left(\frac{N}{x+i},1\right) \right)\nonumber \\
                            &+& \left( F_{n-1}\left(\frac{N}{x+\ell_1},y_1\right) - F_{n-1}\left(\frac{N}{x^*+\ell_1},y_1\right) \right). \label{2.8'}
\end{eqnarray}
Now for each $B \in {\mathcal B}^2_{I}$ one has
\begin{equation} \label{2.9}
\frac{N}{(N+1)^2 \log \left(\frac{N+1}{N}\right) } \overline{\lambda}(B) \leq \overline{G_N}(B) \leq \frac{1}{N \log \left(\frac{N+1}{N}\right) } \overline{\lambda}(B).
\end{equation}
Now from (\ref{1.7'}), (\ref{2.9}) and the fact that $\overline{T_N}$ is $\overline{G_N}$-invariant, it follows that:
\begin{eqnarray}
\sum_{i \geq \ell_1} \left( \frac{N}{x^*+i} - \frac{N}{x+i} \right) &=&
N \sum_{i \geq \ell_1} \overline{\lambda} \left( \left(\left[0, \frac{N}{x+i}\right], \left[0, \frac{N}{x^*+i}\right] \right) \times [0,1]\right) \nonumber \\
&\leq& (N+1)^2 \log \left(\frac{N+1}{N}\right) \sum_{i \geq \ell_1} \overline{G_N} \left( (x^*,x) \times I_N(i) \right) \nonumber\\
&\leq& \frac{(N+1)^2}{N} (x-x^*) \sum_{i \geq \ell_1} {\lambda} \left(I_N(i)\right) \nonumber \\
&=& \frac{(N+1)^2}{N} (x-x^*) \frac{N}{\ell_1} \leq \frac{(N+1)^3}{N} (x-x^*) y \nonumber \\
&=& \frac{(N+1)^3}{N} \overline{\lambda} \left( \Delta_{x,y} \setminus \Delta_{x^*,y}  \right).  \label{2.10}
\end{eqnarray}
For every $2 \leq k \leq d$, a similar analysis leads to
\[
\sum_{i \geq \ell_k} \left| [i, \ell_{k-1}, \ldots, x^*+\ell_{1}]_N - [i, \ell_{k-1}, \ldots, x+\ell_{1}]_N \right| \leq \qquad \qquad\qquad \qquad\qquad
\]
\[
(N+1)^2 \log \left(\frac{N+1}{N}\right) \sum_{i \geq \ell_k} \overline{G_N} \left( \left(\overline{T_N}\right)^k \left( [i, \ell_{k-1}, \ldots, x^*+\ell_{1}]_N, \right. \right.
%\]
%
%\[
\left.\left.[i, \ell_{k-1}, \ldots, x+\ell_{1}]_N\right) \times [0,1] \right)
\]
\[
\leq \frac{(N+1)^2}{N} \sum_{i \geq \ell_k} \overline{\lambda} \left( (x^*,x) \times I_N\left(\ell_1, \ldots, \ell_{k-1},i\right) \right) \ \
\]
\[
\leq \frac{(N+1)^2}{N} (x-x^*) \lambda \left(I_N\left(\ell_1, \ldots, \ell_{k-1}\right)\right) \qquad \qquad \
\]
\begin{equation}\label{2.11}
\leq \frac{(N+1)^3}{N} (x-x^*) y = \frac{(N+1)^3}{N} \overline{\lambda} \left( \Delta_{x,y} \setminus \Delta_{x^*,y}  \right).
\end{equation}
Since $F_n(x,1) = F_n(x)$, from Theorem \ref{G-K1} it follows that
\begin{eqnarray}
\sum_{i \geq \ell_1} \left( F_{n-1}\left( \frac{N}{x^*+i},1\right) - F_{n-1}\left( \frac{N}{x+i},1\right) \right) = \qquad \qquad \qquad\qquad\qquad \nonumber \\
\sum_{i \geq \ell_1} \left( G_N \left( \left[\frac{N}{x+i},\frac{N}{x^*+i} \right]\right) +
\left( \frac{N}{x^*+i} - \frac{N}{x+i}\right){\mathcal O} \left(q^{n-1}\right) \right) = \qquad \ \nonumber \\
\frac{1}{\log \left(\frac{N+1}{N}\right)} \sum_{i \geq \ell_1} \log \left( \frac{x^*+i+1}{x+i+1} \cdot \frac{x+i}{x^*+i} \right) + \qquad\qquad\qquad\qquad\qquad\quad \ \nonumber \\
+ \frac{(N+1)^3}{N} \overline{\lambda} \left( \Delta_{x,y} \setminus \Delta_{x^*,y}  \right) {\mathcal O} \left(q^{n-1}\right) \leq \qquad\qquad\qquad\qquad\qquad\quad\quad \quad\nonumber \\
\frac{1}{\log \left(\frac{N+1}{N}\right)} \cdot \log \left(\frac{x+\ell_1}{x^*+\ell_1}\right) +
\frac{(N+1)^3}{N} \overline{\lambda} \left( \Delta_{x,y} \setminus \Delta_{x^*,y}  \right) {\mathcal O} \left(q^{n-1}\right). \qquad \label{2.12}
\end{eqnarray}
Now from (\ref{2.8'}), (\ref{2.12}), we have:
\begin{eqnarray*}
&& F_{n}(x,y) - F_{n}(x^*,y)=
\frac{1}{\log \left(\frac{N+1}{N}\right)}  \log \left(\frac{x+\ell_1}{x^*+\ell_1}\right) + \\
&& \frac{(N+1)^3}{N} \overline{\lambda} \left( \Delta_{x,y} \setminus \Delta_{x^*,y}  \right) {\mathcal O} \left(q^{n-1}\right)+\\
&& \frac{1}{\log \left(\frac{N+1}{N}\right)} \log \left(\frac{\ell_2+\frac{N}{x+\ell_1}}{\ell_2+\frac{N}{x+\ell_1}}\right)+
\frac{(N+1)^3}{N} \overline{\lambda} \left( \Delta_{x,y} \setminus \Delta_{x^*,y}  \right) {\mathcal O} \left(q^{n-2}\right) +\\
&& F_{n-2}\left(\frac{N}{\ell_2+\frac{N}{\ell_1+x}}, y_2\right) - F_{n-2}\left(\frac{N}{\ell_2+\frac{N}{\ell_1+x^*}}, y_2\right).
\end{eqnarray*}
Applying (\ref{2.8'}) $d$-times and tacking into account that $y_d = 0$, we get
\begin{eqnarray*}
F_{n}(x,y) - F_{n}(x^*,y)= \nonumber \\
\frac{1}{\log \left(\frac{N+1}{N}\right)} \log  \left( \frac{x+\ell_1}{x^*+\ell_1} \cdot \frac{[x+\ell_1]_N+\ell_2}{[x^*+\ell_1]_N+\ell_2} \cdot \frac{[\ell_{d-1}, \ldots, \ell_2, x+\ell_1]_N+\ell_d}{[\ell_{d-1}, \ldots, \ell_2, x^*+\ell_1]_N+\ell_d} \right) + \nonumber \\
\frac{(N+1)^3}{N} \overline{\lambda} \left( \Delta_{x,y} \setminus \Delta_{x^*,y}  \right) \left({\mathcal O} \left(q^{n-1}\right) + \ldots+{\mathcal O} \left(q^{n-d}\right) \right).
\end{eqnarray*}
If $p_d$ and $q_d$ are as in (\ref{1.4}) and (\ref{1.5}) with $a_1=x+\ell_1$ and $a_d=\ell_d$, $d \geq 2$, then
$[\ell_i, \ldots, \ell_2, x+\ell_1]_N =N \displaystyle\frac{q_{i-1}}{q_i}$, $i=1,\ldots,d$.
Since $q_0=1$, thus we have
\[
(x+\ell_1)([x+\ell_1]_N+\ell_2)\cdots([\ell_{d-1}, \ldots, \ell_2, x+\ell_1]_N+\ell_d)=q_d.
\]
Let $p^*_d$ and $q^*_d$ are as in (\ref{1.4}) and (\ref{1.5}), with $a_1=x^*+\ell_1$ and $a_d=\ell_d$, $d \geq 2$. Note that $p_d = p^*_d$.
Thus we find that
\begin{eqnarray*}
\frac{(x+\ell_1)([x+\ell_1]_N+\ell_2)\cdots([\ell_{d-1}, \ldots, \ell_2, x+\ell_1]_N+\ell_d)}
{(x^*+\ell_1)([x^*+\ell_1]_N+\ell_2)\cdots([\ell_{d-1}, \ldots, \ell_2, x^*+\ell_1]_N+\ell_d)} = \frac{q_d}{q^*_d} = \nonumber \\
\frac{p^*_d}{q^*_d} \frac{q_d}{p_d} = \frac{x+\ell_1+[\ell_2, \ldots, \ell_d]_N}{x^*+\ell_1+[\ell_2, \ldots, \ell_d]_N}=
\frac{x+\frac{N}{y}}{x^*+\frac{N}{y}} = \frac{xy+N}{x^*y+N}. \qquad\qquad\quad
\end{eqnarray*}
Therefore,
\begin{eqnarray*}
F_{n}(x,y) - F_{n}(x^*,y)= \frac{1}{\log \left(\frac{N+1}{N}\right)} \log \left( \frac{xy+N}{x^*y+N} \right) \nonumber \\
+ \frac{(N+1)^3}{N} \overline{\lambda} \left( \Delta_{x,y} \setminus \Delta_{x^*,y}  \right) {\mathcal O} \left(q^{n-d}\right),
\end{eqnarray*}
which completes the proof.
\hfill $\Box$

\noindent \textbf{Proof of Theorem \ref{G-K2}.}
Since $\Delta_{x,p_d/q_d} \subset \Delta_{x,y}$ and $F_n(x,y)=\overline{\lambda} \left(\left(\overline{T_N}\right)^{-n}(\Delta_{x,y})\right)$, from (\ref{1.6}), (\ref{2.9}) and the fact that $\overline{T_N}$ is $\overline{G_N}$-invariant, we find that
\begin{eqnarray} \label{2.16}
F_{n}(x,y) - F_{n}(x,\frac{p_d}{q_d}) = \overline{\lambda}\left(\left(\overline{T_N}\right)^{-n}(\Delta_{x,y}) \setminus \left(\overline{T_N}\right)^{-n}(\Delta_{x,p_d/q_d}) \right) \qquad  \nonumber \\
\leq \frac{(N+1)^2}{N} \log \left( \frac{N+1}{N} \right) \overline{G_N}\left(\left(\overline{T_N}\right)^{-n}(\Delta_{x,y}) \setminus \left(\overline{T_N}\right)^{-n}(\Delta_{x,p_d/q_d}) \right) \nonumber \\
\leq \left( \frac{N+1}{N} \right)^2 \overline{\lambda} \left( [0,x] \times \left[ \frac{p_d}{q_d}, y\right] \right) \qquad\qquad\qquad\qquad \qquad\qquad\qquad \  \nonumber \\
\leq \left( \frac{N+1}{N} \right)^2 x \left| y - \frac{p_d}{q_d} \right| \leq \left( \frac{N+1}{N} \right)^2 \cdot \frac{x}{N^n}.  \qquad\qquad\qquad\qquad\qquad
\end{eqnarray}
Since for every fixed $x \in [0,1]$ the function $y \mapsto \log \left(\frac{xy+N}{N}\right)$ is a differentiable on $[0,1]$,
by the \textit{Mean Value Theorem} we have
\begin{eqnarray} \label{2.17}
\left| \log \left(\frac{xy+N}{N}\right) - \log \left(\frac{x\frac{p_d}{q_d}+N}{N}\right) \right| &=& \left| y-\frac{p_d}{q_d} \right| \cdot
\left| \frac{x}{x \xi +N} \right| \nonumber \\
&\leq& x \cdot \left| y - \frac{p_d}{q_d} \right| \leq \frac{x}{N^n},
\end{eqnarray}
where $p_d/q_d \leq \xi \leq y$.
From Lemma \ref{lema2.3}, (\ref{2.16}) and (\ref{2.17}), we have
\begin{eqnarray*}
&&\left| F_{n}(x,y) - \frac{1}{\log\left(\frac{N+1}{N}\right)} \log \left( \frac{xy+N}{N}\right) \right|
\leq \left| F_{n}(x,y) - F_{n}\left(x,\frac{p_d}{q_d}\right)\right|  \\
&&+ \left| F_n(x,\frac{p_d}{q_d}) - F_n(0,\frac{p_d}{q_d}) - \frac{1}{\log\left(\frac{N+1}{N}\right)} \log \left( \frac{x\frac{p_d}{q_d}+N}{N} \right) \right| \\
&&+ \frac{1}{\log\left(\frac{N+1}{N}\right)} \left| \log \left(\frac{xy+N}{N}\right) - \log \left(\frac{x\frac{p_d}{q_d}+N}{N}\right) \right|  \\
&&\leq  \left( \frac{N+1}{N} \right)^2 \cdot \frac{x}{N^n} + C q^{n-d} + \frac{x}{N^n} 
\end{eqnarray*}
which completes the proof.
\hfill $\Box$
\section{Improving result}

In this section we shall estimate the error term
\begin{equation*}
e_{n,a} (x,y) = G_{N,a} \left( T^n_N \in [0,x], s_{n,a} \in [0,y] \right) - \frac{1}{\log\left(\frac{N+1}{N}\right)} \log\left( \frac{xy+N}{N} \right)
\end{equation*}
for any $a, x, y \in I$ and $n \in \mathbb{N}$.

In the main result of this section, Theorem \ref{th.4.7}, we shall derive lower and upper bounds (not depending on $a \in I$) of the supremum
\begin{equation}\label{4.2}
\sup_{x \in I, y \in I} |e_{n,a} (x, y)|, \quad a \in I,
\end{equation}
which provide a more refined estimate of the convergence rate involved.
First, we obtain a lower bound for the error, which suggests the exact convergence rate of
$G_{N,a} \left( s_{n,a} \in [0,y] \right)$ to $G_{N} \left([0,y] \right)$ as $n \rightarrow \infty$ for all $a \in I$.

\begin{theorem} \label{Th.4.1}
Whatever $a \in I$ we have
\begin{equation}\label{4.3}
\frac{1}{2} v_{N(n)}(1) \leq \sup_{y \in I} \left|G_{N,a} \left( s_{n,a} \in [0,y] \right) - G_{N} \left([0,y] \right) \right|
\end{equation}
for all $n \in \mathbb{N}_+$.
\end{theorem}
\noindent \textbf{Proof.}
First, using the continuity of $G_{N} \left([0,y] \right)$ with respect to $y$ it is easy to see that
\begin{equation*}
\sup_{y \in I} \left|G_{N,a} \left( s_{n,a} \in [0,y] \right) - G_{N} \left([0,y] \right) \right| =
\sup_{y \in I} \left|G_{N,a} \left( s_{n,a} \in [0,y) \right) - G_{N} \left([0,y] \right) \right|
\end{equation*}
for all $a \in I$ and $n \in \mathbb{N}$. Second, whatever $s \in I$ we have
\begin{eqnarray*}
  G_{N,a} (s_{n,a}=s) &=& G_{N,a} \left( s_{n,a} \in [0,s] \right) -  G_{N} \left([0,s] \right) \nonumber \\
  &-& \left(G_{N,a} \left( s_{n,a} \in [0,s) \right) - G_{N} \left([0,s] \right) \right) \nonumber \\
  &\leq& \sup_{y \in I} \left|G_{N,a} \left( s_{n,a} \in [0,y] \right) - G_{N} \left([0,y] \right) \right| \nonumber \\
  &+& \sup_{y \in I} \left|G_{N,a} \left( s_{n,a} \in [0,y) \right) - G_{N} \left([0,y] \right) \right| \nonumber \\
  &=& 2 \sup_{y \in I} \left|G_{N,a} \left( s_{n,a} \in [0,y] \right) - G_{N} \left( [0,y] \right) \right|.
\end{eqnarray*}
Hence
\begin{equation*}
\sup_{y \in I} \left|G_{N,a} \left( s_{n,a} \in [0,y] \right) - G_{N} \left( [0,y] \right) \right| \geq
\frac{1}{2} \sup_{s \in I} G_{N,a} \left( s_{n,a} = s \right),
\end{equation*}
for all $a \in I$ and $n \in \mathbb{N}$. For any $n \in \mathbb{N}_+$ we have
\begin{equation}\label{4.7}
U^n f(y) = \sum_{i_1,\ldots,i_n \in \Lambda} v_{i_1\ldots i_n} (x) f\left( u_{i_n \ldots i_1} (x)\right)
\end{equation}
where $u_{i_n \ldots i_1} = u_{i_n} \circ \ldots \circ u_{i_1}$, $v_{i_1\ldots i_n} (x)= v_{i_1}(x) v_{i_2}(u_{i_1}(x)) \ldots v_{i_n} (u_{i_{n-1} \ldots i_1} (x))$, $n \geq 2$, and the functions $u_i$ and $v_i$, $i \in \Lambda$, are defined in (\ref{1.17'}) and (\ref{1.17''}).
Now, we obtain
\begin{equation*}
G_{N,a} \left( s_{n,a} = [i_n, \ldots, i_2, i_1+a]_N \right) = G_{N,a} \left( I_N \left( i^{(n)}\right) \right) =  v_{i_1\ldots i_n}(a), \ n\geq 2,
\end{equation*}
\begin{equation*}
G_{N,a} \left( s_{1,a} = \frac{N}{i_1 +a} \right) = G_{N,a} \left( I_N \left( i_{1}\right) \right) = v_{i_1}(a)
\end{equation*}
for all $a \in I$ and $i_1,\ldots,i_n \in \Lambda$.

Putting
\begin{equation*}
\frac{p_n(i_1, \ldots, i_n)}{q_n(i_1, \ldots, i_n)} = [i_1, \ldots, i_n]_N
\end{equation*}
we get
\begin{eqnarray} \label{4.11}
v_{i_1\ldots i_n} (x) &=& \frac{(x+N)N^{n-1}}{q_{n-1}(i_2, \ldots, i_n)(x+i_1) + p_{n-1}(i_2, \ldots, i_n)} \nonumber \\
                      &\times& \frac{1}{q_{n}(i_2, \ldots, i_n, N)(x+i_1) + p_{n}(i_2, \ldots, i_n, N)}
\end{eqnarray}
for all $i_n \in \Lambda$, $n \geq 2$, and $a \in I$.

By (\ref{4.11}) we have
\begin{equation}\label{4.12}
\sup_{s \in I} G_{N,a} \left( s_{n,a} = s \right) = v_{N(n)}(a)
\end{equation}
where we write $N(n)$ for $(i_1, \ldots, i_n)$ with $i_1=\ldots=i_n=N$, $n \in \mathbb{N}_+$.

Also by (\ref{4.11}) we have
\begin{eqnarray*}
v_{N(n)} (a) &=& \frac{(a+N)N^{n-1}}{q_{n-1}(\underbrace{N, \ldots, N}_{(n-1) \ times})(a+N) + p_{n-1}(\underbrace{N, \ldots, N}_{(n-1) \ times})} \nonumber \\
                      &\times& \frac{1}{q_{n}(\underbrace{N, \ldots, N, N}_{n \ times})(a+N) + p_{n}(\underbrace{N, \ldots, N, N}_{n \ times})}.
\end{eqnarray*}
It is easy to see that $v_{N(n)} (\cdot)$ is a decreasing function. Therefore
\begin{equation*}
\sup_{s \in I} G_{N,a} \left( s_{n,a} = s \right) \geq v_{N(n)}(1)
\end{equation*}
for all $a \in I$.
\hfill $\Box$

\begin{theorem} \label{th.4.2}
(The lower bound)
Whatever $a \in I$ we have
\begin{equation*}
\frac{1}{2} v_{N(n)}(1) \leq \sup_{x \in I, y \in I} \left|G_{N,a} \left( T^n_N \in [0,x], s_{n,a} \in [0,y] \right) - \frac{1}{\log\left(\frac{N+1}{N}\right)} \log\left( \frac{xy+N}{N} \right) \right|
\end{equation*}
for all $n \in \mathbb{N}_+$.
\end{theorem}
\noindent \textbf{Proof.} Whatever $a \in I$ and $n \in \mathbb{N}_+$, by Theorem \ref{Th.4.1} we have
\begin{eqnarray*}
&\sup_{x \in I, y \in I}& \left|G_{N,a} \left( T^n_N \in [0,x], s_{n,a} \in [0,y] \right) - \frac{1}{\log\left(\frac{N+1}{N}\right)} \log\left( \frac{xy+N}{N} \right) \right| \nonumber \\
&\geq& \sup_{y \in I} \left|G_{N,a} \left( T^n_N \in I, s_{n,a} \in [0,y] \right) - \frac{1}{\log\left(\frac{N+1}{N}\right)} \log\left( \frac{y+N}{N} \right) \right| \nonumber \\
&=& \sup_{y \in I} \left|G_{N,a} \left( s_{n,a} \in [0,y] \right) - G_{N} \left( [0,y] \right) \right| \geq \frac{1}{2} v_{N(n)}(1).
\end{eqnarray*}
\hfill $\Box$

\begin{remark}
Here
\begin{equation*}
v_{N(n)}(1) = \frac{(1+N)N^{n+1}}{q_{n+1} q_{n+2}}, \ n \in \mathbb{N}_+,
\end{equation*}
where $q_n=N(q_{n-1} + q_{n-2})$, $n \in \mathbb{N}_+$, with $q_{-1} = 0$ and $q_0 = 1$.
It is easy to see that
\begin{equation*}
q_n = \frac{1}{\sqrt{N^2 + 4N}} \left[ \left( \frac{N+\sqrt{N^2 + 4N}}{2} \right)^{n+1} - \left( \frac{N-\sqrt{N^2 + 4N}}{2} \right)^{n+1} \right].
\end{equation*}
It should be noted that Theorem \ref{th.4.2} in connection with the limit
\begin{equation*}
\lim_{n \rightarrow \infty} \left( \frac{1}{2} v_{N(n)}(1) \right)^{1/n} = \frac{2}{N+\sqrt{N^2 + 4N}+2}
\end{equation*}
leads to an estimate of the order of magnitude of the error $e_{n,a} (x,y)$.

It is known that for the RCF-expansion \cite{I-1997} the exact order of magnitude of the supremum there is $\mathcal{O}(g^{2n})$ with $g = (\sqrt{5}-1)/2$, $g^2 = (3-\sqrt{5})/2 = 0.38196\ldots$. Note that for $N=1$,
$\displaystyle\lim_{n \rightarrow \infty} \left( \frac{1}{2} v_{N(n)}(1) \right)^{1/n} = g^2$.
\end{remark}
In what follows we study the transition operator associated with the RSCC underlying N-continued fraction on the Banach space of complex-valued functions of bounded variation. The characteristic properties of this operator are used to derive an explicit upper bound for $\displaystyle\sup_{x \in I, y \in I} |e_{n,a}(x,y)|$, $a \in I$.

Let $BV(I)$ the Banach space of complex-valued functions $f$ of bounded variation on $I$ under the norm
\begin{equation*}
\left\| f \right\|_\mathrm{v} := \mathrm{var}f + |f|.
\end{equation*}
Remember that the variation ${\rm var}_{A}f$ over $A \subset I$ of $f \in B(I)$ is defined as
\begin{equation*}
\sup \sum^{k}_{i=1} |f(t_{i}) - f(t_{i+1})|
\end{equation*}
the supremum being taken over all $t_1 < \cdots < t_k \in A$, $k \geq 2$.
We write simply $\mathrm{var} f$ for $\mathrm{var}_I f$ and if $\mathrm{var} f < \infty$, then $f$ is called a function of bounded variation.

We start by proving the following elementary result.
\begin{proposition} \label{prop4.4}
For any $f \in BV(I)$ we have
\begin{equation*}
\mathrm{var}\,Uf \leq \frac{1}{N+1} \cdot \mathrm{var } f.
\end{equation*}
\end{proposition}
\noindent \textbf{Proof.}
Put $v_i(x) = (x+N)\left( \displaystyle \frac{1}{x+i} - \displaystyle \frac{1}{x+i+1}\right)$, $i \in \Lambda$.
We have
\[
v'_i(x) = \frac{1}{(x+i)(x+i+1)} \left[ 1 - \frac{(x+N)(2x+2i+1)}{(x+i)(x+i+1)}\right], \, i \in \Lambda.
\]
It follows that $v'_i(x) < 0$, $x \in I$, $i \in \Lambda$. Hence
\begin{eqnarray*}
  \mathrm{var} \, v_i &=& v_i(0) - v_i(1), i \in \Lambda \\
  |v_i| &=& \sup_{x \in I} v_i(x) = v_i(0) = \frac{N}{i(i+1)}, i \in \Lambda.
\end{eqnarray*}
Thus
\begin{equation*}
\sup_{i \in \Lambda} |v_i| = v_N (0) = \frac{1}{N+1}.
\end{equation*}
Also,
\begin{equation*}
\sum_{i \in \Lambda} \mathrm{var} \, v_i = \sum_{i \in \Lambda} \left( \frac{N}{i(i+1)} - \frac{N+1}{(i+1)(i+2)} \right) = 0.
\end{equation*}
We have
\begin{eqnarray*}
\mathrm{var} \, Uf &=& \mathrm{var} \sum_{i \in \Lambda} v_i \cdot (f \circ u_i) \leq \sum_{i \in \Lambda} \mathrm{var} \left( v_i \cdot (f \circ u_i)\right) \\
&\leq& \sum_{i \in \Lambda} |v_i| \mathrm{var} (f \circ u_i) +  \sum_{i \in \Lambda} |f \circ u_i| \mathrm{var} \, v_i \\
&\leq& \left( \sup_{i \in \Lambda} |v_i| \right)\sum_{i \in \Lambda} \mathrm{var} (f \circ u_i) + |f| \sum_{i \in \Lambda} \mathrm{var} \, v_i \\
&\leq& \frac{1}{N+1} \cdot \mathrm{var } f
\end{eqnarray*}
because
\begin{equation*}
\sum_{i \in \Lambda} \mathrm{var} (f \circ u_i) = \sum_{i \in \Lambda} \mathrm{var}_ {\left[ \frac{N}{i+1}, \frac{N}{i} \right]} \, f = \mathrm{var} \, f.
\end{equation*}
\hfill $\Box$
\begin{corollary} \label{cor.4.5}
For any $f \in BV(I)$ and for all $n \in \mathbb{N}$ we have
\begin{eqnarray}
 \mathrm{var}\,U^n f &\leq& \frac{1}{(N+1)^n} \cdot \mathrm{var } f, \label{4.19} \\
  \left| U^n f - U^{\infty} f \right| &\leq& \frac{1}{(N+1)^n} \cdot \mathrm{var } f. \label{4.20}
\end{eqnarray}
\end{corollary}
\noindent \textbf{Proof.} Note that for any $f \in BV(I)$ and $u \in I$ we have
\begin{eqnarray*}
|f(u)| -  \left| \int_{I} f(x) G_N (\mathrm{dx}) \right| &\leq& \left| f(u)- \int_{I} f(x) G_N (\mathrm{dx}) \right| \\
&=&  \left| \int_{I}(f(u) - f(x)) G_N (\mathrm{dx})\right| \leq \mathrm{var } f,
\end{eqnarray*}
whence
\begin{equation}\label{4.21}
|f| \leq  \left| \int_{I} f(x) G_N (\mathrm{dx}) \right| + \mathrm{var } f, \quad f \in BV(I).
\end{equation}
Finally, (\ref{1.18}) and (\ref{4.21}) imply that
\begin{equation*}
\left| U^n f - U^{\infty} f \right| \leq  \mathrm{var } \left( U^n f - U^{\infty} f \right) = \mathrm{var } \, U^n f.
\end{equation*}
for all $n \in \mathbb{N}$ and $f \in BV(I)$, which leads to (\ref{4.20}).
\hfill $\Box$
\begin{theorem} \label{th.4.6}
(The upper bound)
Whatever $a \in I$ we have
\begin{equation*}
\sup_{x \in I, y \in I} \left|G_{N,a} \left( T^n_N \in [0,x], s_{n,a} \in [0,y] \right) - \frac{1}{\log\left(\frac{N+1}{N}\right)} \log\left( \frac{xy+N}{N} \right) \right| \leq \frac{1}{(N+1)^n}
\end{equation*}
for all $n \in \mathbb{N}$.
\end{theorem}
\noindent \textbf{Proof.} Let $F_{n,a}(y) = G_{N,a} (s_{n,a} \leq y)$ and $H_{n,a} (y) = F_{n,a}(y) - G_N ([0,y])$, $a,y \in I$, $n \in \mathbb{N}$. As we have noted $U$ is the transition operator of the Markov chain $(s_{n,a})_{n \in \mathbb{N}}$. For any $y \in I$ consider the function $f_y$ defined on $I$ as
\begin{equation*}
f_{y}(a):=
\left\{
\begin{array}{ll}
{1} & { \mbox{if } \, 0 \leq a \leq y,}\\
{0} & \mbox{if } \, y < a \leq 1.
\end{array}
\right.
\end{equation*}
Hence
\begin{equation*}
U^n f_y (a) = E_a\left( \left. f_y(s_{n,a})\right| s_{0,a} = a \right) = G_{N,a} (s_{n,a} \leq y)
\end{equation*}
for all $a,y \in I$, $n \in \mathbb{N}$.
As
\begin{equation*}
U^{\infty} f_y = \int_{I} f_y(a) G_N (\mathrm{da}) = G_N ([0,y]), \quad y \in I.
\end{equation*}
It follows from Corollary \ref{cor.4.5} that
\begin{eqnarray} \label{4.24}
  |H_{n,a}(y)| &=& \left| G_{N,a} (s_{n,a} \leq y) - G_N ([0,y]) \right| \nonumber \\
               &=& \left| U^n f_y (a) - U^{\infty} f_y \right| \leq \frac{1}{(N+1)^n} \mathrm{var } \, f_y = \frac{1}{(N+1)^n}
\end{eqnarray}
for all $a,y \in I$, $n \in \mathbb{N}$.
By (\ref{1.19}), for all $a \in I$, $x,y \in I$ and $n \in \mathbb{N}$ we have
\begin{eqnarray*}
  G_{N,a} \left( T^n_N \in [0,x], s_{n,a} \in [0,y] \right) = \int^{y}_{0}  G_{N,a} \left( \left.T^n_N \in [0,x] \right| s_{n,a} = z \right)
  \mathrm{d F_{n,a}(z)}  \nonumber \\
   = \int^{y}_{0} \frac{(z+N)x}{zx+N} \mathrm{d F_{n,a}(z)} = \int^{y}_{0} \frac{(z+N)x}{zx+N} \mathrm{d G_{N}(z)} + \int^{y}_{0} \frac{(z+N)x}{zx+N} \mathrm{d H_{n,a}(z)} \nonumber \\
   = \frac{1}{\log\left(\frac{N+1}{N}\right)} \log\left( \frac{xy+N}{N} \right) + \frac{(z+N)x}{zx+N} \left.H_{n,a}(z)\right|^{y}_{0}
   - \int^{y}_{0} \frac{Nx(1-x)}{(zx+N)^2} H_{n,a}(z)\mathrm{dz}.
\end{eqnarray*}
Hence, by (\ref{4.24})
\begin{eqnarray*}
\left|G_{N,a} \left( T^n_N \in [0,x], s_{n,a} \in [0,y] \right) - \frac{1}{\log\left(\frac{N+1}{N}\right)} \log\left( \frac{xy+N}{N} \right) \right|\\
\leq \frac{1}{(N+1)^n} \left( \frac{(y+N)x}{xy+N} - \frac{(1-x)xy}{xy+N} \right) = \frac{x}{(N+1)^n} \leq \frac{1}{(N+1)^n}
\end{eqnarray*}
for all $a,x,y \in I$ and $n \in \mathbb{N}$.
\hfill $\Box$

Combining Theorem \ref{th.4.2} with Theorem \ref{th.4.6} we obtain Theorem \ref{th.4.7}.
\begin{theorem} \label{th.4.7}
Whatever $a \in I$ we have
\begin{eqnarray*}
\frac{1}{2} v_{N(n)}(1) \leq \qquad  \qquad\qquad\qquad\qquad\qquad\qquad\qquad\qquad\qquad\qquad\qquad\qquad\qquad\\
 \sup_{x \in I, y \in I} \left|G_{N,a} \left( T^n_N \in [0,x], s_{n,a} \in [0,y] \right) - \frac{1}{\log\left(\frac{N+1}{N}\right)} \log\left( \frac{xy+N}{N} \right) \right| \leq \frac{1}{(N+1)^n}
\end{eqnarray*}
for all $n \in \mathbb{N}_+$.
\end{theorem}
\begin{remark}
Theorem \ref{th.4.7} implies that the convergence rate is $\mathcal{O}(\alpha^n)$, with
\begin{equation*}
\frac{2}{N+\sqrt{N^2+4N}+2} \leq \alpha \leq \frac{1}{N+1}.
\end{equation*}
\end{remark}
For example, we have

\begin{center}
\begin{tabular}{|l|l|}
  \hline
  $N=1$ & $g^2=0.381966 \leq \alpha \leq 0.5$\\\hline
  $N=2$ & $0.267949192 \leq \alpha \leq 0.333333\ldots$\\\hline
  $N=5$ & $0.145898033 \leq \alpha \leq 0.166666\ldots$\\\hline
  $N=10$ & $0.083920216 \leq \alpha \leq 0.090909\ldots$\\\hline
  $N=100$ & $0.009804864 \leq \alpha \leq 0.00990099$\\\hline
  $N=1000$ & $0.000998004 \leq \alpha \leq 0.000999$\\\hline
  $N=10000$ & $0.00009998 \leq \alpha \leq 0.00009999$\\
  \hline
\end{tabular}
\end{center}


\begin{thebibliography}{[01]}
%
\bibitem{Burger-2008} Burger, E. B., Gell-Redman, J., Kravitz, R., Walton, D. and Yates, N., \textit{Shrinking the period lengths of continued fractions while still capturing convergents.} J. Number Theory, \textbf{128} (2008) 1,  144-153.
%

\bibitem{DK-1994} Dajani, K. and Kraaikamp, C., \textit{Generalization of a theorem by Kusmin}, Monatsh. Math., \textbf{118} (1994) 55-73.
%
\bibitem{DK-1999} Dajani, K. and Kraaikamp, C., \textit{A Gauss-Kuzmin theorem for optimal continued fractions.}, Trans. Amer. Math. Soc. \textbf{351} (1999) 2055-2079.
%
\bibitem{DKW-2013} Dajani, K., Kraaikamp, C. and Van der Wekken, N., \textit{Ergodicity of N-continued fraction expansions.} J. Number Theory, \textbf{133} (2013) 9, 3183-3204.
%
\bibitem{I-1997} Iosifescu, M., \textit{On the Gauss-Kuzmin-L\'evy theorem, III.} Rev. Roumaine de Math. Pures Appl. \textbf{42} (1) (1997) 71-88.
%
\bibitem{IG-2009} Iosifescu, M. and Grigorescu, S., \textit{Dependence With Complete Connections and its Applications}, Cambridge Tracts in Mathematics \textbf{96}, Cambridge Univ.Press, Cambridge, 2nd edition, 2009.
%
\bibitem{IK-2002} Iosifescu, M. and Kraaikamp, C., \textit{Metrical Theory of Continued Fractions.} Kluwer Academic Publishers, Dordrecht, 2002.
%
\bibitem{Lascu-2016} Lascu, D., \textit{Dependence with complete connections and the Gauss-Kuzmin theorem for $N$-continued fractions.} J. Math. Anal. Appl. \textbf{444} (2016) 610–623.
%
\bibitem{Lascu-2017} Lascu, D., \textit{Metric properties of $N$-continued fractions.} Math. Reports \textbf{19}(69), 2 (2017) 165-181.
%
\bibitem{Nakada-1981} Nakada, H., \textit{Metrical theory for a class of continued fraction transformations and their natural extensions.} Tokyo J. Math. \textbf{4} (1981), 2, 399-426.
%
\bibitem{Sebe-2000} Sebe, G. I., \textit{A two-dimensional Gauss-Kuzmin theorem for singular continued fractions.} Indag. Mathem., N.S., \textbf{11} (4) (2000) 593-605.
%
\bibitem{Sebe-2001} Sebe, G.I., \textit{On convergence rate in the Gauss-Kuzmin problem for grotesque continued fractions}, Monatsh. Math. \textbf{133} (2001) 241-254.

\end{thebibliography}
\end{document}